\documentclass[a4paper,12pt]{amsart}
\usepackage{mathrsfs}
\usepackage{amsfonts}
\usepackage{amssymb}
\usepackage{ifthen}
\usepackage{graphicx}
\usepackage[usenames]{color}
\usepackage{amsmath}
\usepackage{amscd}
\usepackage{float}
\usepackage{color}

\font\ff=eusm10 scaled 1200

\def\K{\hbox{\ff K}}

\nonstopmode \numberwithin{equation}{section}
\newtheorem{thm}{Theorem} %[section]
\newtheorem{lem}{Lemma} %[section]
\newtheorem{cor}{Corollary} %[section]
\newtheorem{conj}[equation]{Conjecture}

\theoremstyle{definition}
\newtheorem{defn}{Definition}%[section]
\newtheorem{prob}[equation]{Problem}
\newtheorem{ques}[equation]{Question}
\newtheorem{rem}{Remark}%[section]
\newtheorem{exam}[equation]{Example}

\newcounter {own}
\def\theown {\thesection       .\arabic{own}}

\newenvironment{pf}[1][]{%
 \vskip 3mm
 \noindent
 \ifthenelse{\equal{#1}{}}%
  {{\slshape Proof. }}%
  {{\slshape #1.} }%
 }%
{\qed\bigskip}

\newcounter{alphabet}
\newcounter{tmp}

\makeatletter
\newcommand{\Ref}[1]{\@ifundefined{r@#1}{}{\setcounter{tmp}{\ref{#1}}\Alph{tmp}}}
\makeatother
\newcommand{\IR}{{\mathbb R}}

\newcommand{\ID}{{\mathbb D}}

\def\be{\begin{equation}}
\def\ee{\end{equation}}

\newcommand{\bee}{\begin{enumerate}}
\newcommand{\eee}{\end{enumerate}}

\newcommand{\blem}{\begin{lem}}
\newcommand{\elem}{\end{lem}}
\newcommand{\bthm}{\begin{thm}}
\newcommand{\ethm}{\end{thm}}
\newcommand{\bcor}{\begin{cor}}
\newcommand{\ecor}{\end{cor}}
\newcommand{\beg}{\begin{exam}}
\newcommand{\eeg}{\end{exam}}
\newcommand{\begs}{\begin{examples}}
\newcommand{\eegs}{\end{examples}}
\newcommand{\bdefe}{\begin{defn}}
\newcommand{\edefe}{\end{defn}}
\newcommand{\bprob}{\begin{prob}}
\newcommand{\eprob}{\end{prob}}
\newcommand{\bques}{\begin{ques}}
\newcommand{\eques}{\end{ques}}
\newcommand{\bei}{\begin{itemize}}
\newcommand{\eei}{\end{itemize}}
\newcommand{\bcon}{\begin{conj}}
\newcommand{\econ}{\end{conj}}
\newcommand{\bcons}{\begin{conjs}}
\newcommand{\econs}{\end{conjs}}
\newcommand{\bprop}{\begin{propo}}
\newcommand{\eprop}{\end{propo}}
\newcommand{\br}{\begin{rem}}
\newcommand{\er}{\end{rem}}
\newcommand{\brs}{\begin{rems}}
\newcommand{\ers}{\end{rems}}
\newcommand{\bo}{\begin{obser}}
\newcommand{\eo}{\end{obser}}
\newcommand{\bos}{\begin{obsers}}
\newcommand{\eos}{\end{obsers}}
\newcommand{\bpf}{\begin{pf}}
\newcommand{\epf}{\end{pf}}
\newcommand{\ba}{\begin{array}}
\newcommand{\ea}{\end{array}}
\newcommand{\beq}{\begin{eqnarray}}
\newcommand{\beqq}{\begin{eqnarray*}}
\newcommand{\eeq}{\end{eqnarray}}
\newcommand{\eeqq}{\end{eqnarray*}}

\newcommand{\ra}{\rightarrow}

\newcommand{\ds}{\displaystyle}

\begin{document}
\bibliographystyle{amsplain}
\title {Radii of covering disks for locally univalent harmonic mappings}
%\author{S. Ponnusamy}
%\address{S. Ponnusamy, Department of Mathematics,
%Indian Institute of Technology Madras, Chennai-600 036, India.}
%\email{samy@iitm.ac.in}

\author{Sergey Yu. Graf, Saminathan Ponnusamy, and Victor V. Starkov}

\address{S. Yu. Graf,
Tver State University, ul. Zhelyabova 33, Tver, 170000 Russia.
}
\email{sergey.graf@tversu.ru}

\address{S. Ponnusamy,
Indian Statistical Institute (ISI), Chennai Centre, SETS (Society
for Electronic Transactions and Security), MGR Knowledge City, CIT
Campus, Taramani, Chennai 600 113, India.
}
\email{samy@isichennai.res.in, samy@iitm.ac.in}

\address{V. V. Starkov,
Department of Mathematics, University of Petrozavodsk,
ul. Lenina 33, 185910 Petrozavodsk, Russia
}
\email{vstarv@list.ru }

\keywords{Locally univalent harmonic mappings, linear and affine invariant families, convex and close-to-convex functions,
and covering theorems.
%\\
%${}^{\mathbf{*}}$ Corresponding author
%${} ^\dagger$ {\tt This author is on leave from the Department of Mathematics,
%Indian Institute of Technology Madras, Chennai-600 036, India}
}

%\dedicatory{}
\begin{abstract}
For a univalent smooth mapping $f$ of the unit disk $\ID$ of complex plane onto  the manifold $f(\ID)$, let
$d_f(z_0)$ be the radius of the largest univalent disk on the manifold $f(\ID)$ centered at $f(z_0)$ ($|z_0|<1$).
The main aim of the present article is to investigate how the radius $d_h(z_0)$ varies when the analytic function $h$ is replaced  by
a sense-preserving harmonic function $f=h+\overline{g}$.
The main result includes sharp upper and lower bounds for the quotient $d_f(z_0)/d_h(z_0)$, especially, for
a family of locally univalent $Q$-quasiconformal harmonic mappings $f=h+\overline{g}$ on $|z|<1$. In addition, estimate on
the radius of the disk of convexity of functions belonging to certain linear invariant families of locally univalent
$Q$-quasiconformal harmonic mappings of order $\alpha$ is obtained.
\end{abstract}

\thanks{}

\maketitle
\pagestyle{myheadings}
\markboth{S. Yu. Graf, S. Ponnusamy, and V. V. Starkov}{Linear invariant families of locally univalent harmonic mappings}

\section{Introduction and Main results}

Let $\ID =\{z\in\mathbb C:\, |z|<1\}$ be the unit disk, and $h$ be a smooth univalent mapping of the unit disk $\ID$ onto two-dimensional manifold $M$.
For a point $a\in\ID$, we write  $d_h(z)$ as the radius of the largest univalent disk centered at $h(a)$ on the manifold $M$.
Here a univalent disk on $M$ centered at $h(a)$ means that $h$ maps an open subset of $\ID$ containing the point $a$ univalently onto this disk.

The question about lower estimation of $d_h$ for univalent analytic functions first was considered in papers of
Koebe \cite{Ko} and Bieberbach \cite{Bib} in connection with the well known problem of covering disk in the class
${\mathcal S}$. Here ${\mathcal S}$ denotes the classical family of all normalized univalent (analytic) functions in $\ID$
investigated by a number of researchers (see \cite{Du,Golu66,Pomm}). In the class of analytic functions $h$ in $\ID$ with $h'(0)=1$,
the determination of the exact value of the greatest lower bound of all $d_h$ is one of the most important problems in geometric
function theory of one complex variable. For historical discussion of the attempts of various mathematicians to estimate the
lower bound for $d_h(z)$, we refer to \cite{Minda90} and also \cite{CPR-2014,CPW1,CPW1a} for recent developments.

If $\mathcal{LU}$ denotes the family of functions $h$ analytic and locally univalent $(h'(z)\neq 0)$ in $\ID$, then
the classical Schwarz lemma for analytic functions gives the following well-known sharp upper estimate for $d_h(z)$:
$$d_h(z)\le|h'(z)|(1-|z|^2).
$$
Often the right hand side quantity, namely, $r(h(z), h(\ID))=|h'(z)|(1-|z|^2)$ is referred to as the conformal radius
of the domain $h(\ID)$ at $h(z)$. Sharp and nontrivial lower estimate for $d_h(z)$ was obtained by Pommerenke \cite{Pom64}
in a detailed analysis of what is called \textit{linear invariant families} of locally univalent analytic functions in $\ID$.
Throughout we denote by ${\rm Aut}\,(\ID)$, the set of all conformal automorphisms (M\"obius self-mappings)
$\phi(z)=e^{i\theta} \frac{z+a}{1+\overline{a}z}$, where $|a| <1$ and $\theta\in\IR$, of the unit disk $\ID$.

\bdefe (cf. \cite{Pom64})\label{def1}
A non-empty collection ${\mathfrak M}$ of functions from $\mathcal{LU}$ is called a linear invariant family {\rm (LIF)}
if for each $h\in {\mathfrak M}$, normalized such that $h(z)=z+\sum_{k=2}^\infty a_k(h) z^k$,
%$$
%\aligned
%&h(z)=z+\sum_{k=1}^\infty a_k(h) z^k,\;\;|a_2(h)| \le\alpha,\\
%&\frac{h(\omega(z))-h(\omega(0))}{h'(\omega(0))\omega'(0)} \in{\mathfrak M}\;\;\;\forall\; \omega(z)=e^{i\theta} \frac{z+z_0}{1+\overline{z_0}z}, |z_0| <1, \theta\in\mathbb R.
%\endaligned
%$$
the functions $H_\phi (z)$ defined by
$$ H_\phi (z)=\frac{h(\phi(z))-h(\phi(0))}{h'(\phi(0))\phi'(0)}=z+\ldots ,
$$
belong to ${\mathfrak M}$ for each $\phi\in {\rm Aut}\,(\ID)$.
\edefe

The order of the family ${\mathfrak M}$ is defined to be $\alpha:=\texttt{ord}\,\mathfrak{M}=\sup_{h\in\mathfrak{M}}|a_2(h)|$.
The universal LIF, denoted by ${\mathcal U}_\alpha$, is defined to be the collection of all linear invariant families ${\mathfrak M}$
with order less than or equal to $\alpha$ (see \cite{Pom64}). An interesting fact about the order of a LIF family is that many properties of
it depend only on the order of the family.
It is well-known \cite{Pom64} that ${\mathcal U}_\alpha\ne\emptyset $ if and only if $\alpha\ge 1$. The family
${\mathcal U}_1$ is precisely the family ${\mathcal K}$ of all normalized convex univalent (analytic) functions
whereas ${\mathcal S}\subset {\mathcal U}_2$.

 Note that ${\mathcal U}_\alpha$ is the largest LIF of functions $h$ with the restriction of growth (see \cite{Star90}):
$$|h'(z)|\le\frac{(1+|z|)^{\alpha-1}}{(1-|z|)^{\alpha+1}}.
$$

In \cite{Pom64}, Pommerenke has proved that for each $h\in {\mathcal U}_\alpha$ the following sharp lower estimate of $d_h(z)$ holds:
$$d_h(z)\ge\frac{1}{2\alpha}|h'(z)|(1-|z|^2).
$$

In the present paper we obtain estimate of the functional $d_f(z)$ when instead of analytic functions $h(z)$ we consider harmonic locally
univalent mappings
\be\label{eq1}
f(z)=h(z)+\overline{g(z)}=\sum_{k=1}^\infty\left (a_kz^k+a_{-k}\overline{z}^k\right ),
%\eqno (1)
\ee
i.e. when $\overline{g(z)}$ is added to the functions $h$. In the above decomposition of $f$, the functions $h$ and $g$ are
called the analytic and co-analytic parts of $f$, respectively. We say that a harmonic functions $f=h+\overline{g}$ is sense-preserving if
the Jacobian $J_f(z)=|h'(z)|^2 -|g'(z)|^2$ of $f$ is positive. Lewy's theorem \cite{Le} (see also for
example \cite[Chapter 2, p. 20]{Duren:Harmonic} and \cite{PonRasi2013}) implies that every harmonic function $f$ on $\ID$ is locally one-to-one
and sense-preserving on $\ID$ if and only if $J_f(z)>0$ in $\ID$. Note that $J_f(z)> 0$ in $\ID$ if and only if $h'(z)\neq0$ and there exists
an analytic function $\omega_f$ in $\ID$ such that
\be\label{li1-eq0}
|\omega_f (z)|<1 ~\mbox{ for}~z\in \ID,
\ee
where $\omega_f (z)=g'(z)/h'(z)$. Here $\omega_f$ is referred to as the (complex) \textit{dilatation} of the harmonic mapping $f=h+\overline{g}$.
When it is convenient, we simply use the notation $\omega$ instead of $\omega_f$.

There are different generalizations of the notion of the linear invariant family to the case of harmonic mappings.
For example, the question about a lower estimate of the radius $d_f(0)$ of the univalent disk centered at the origin was
examined by Sheil-Small \cite{Sheil90} in the linear and affine invariant families of univalent harmonic functions $f$.
There are a number of articles in the literature proving such inequalities or studying the related mappings in various settings.
For example, see \cite{CPR-2014,CPW-2011,CPW1,CPW1a,CPW1b,CPW-2013a,CPW-2013b,Star93,S2}, and also the work from \cite{HG} in which one can
obtain a lower bound on the radius for quasi-regular mappings. The concept of linear and affine invariance was also
discussed by Schaubroeck \cite{Schau2000} for the case of locally univalent harmonic mappings.

\bdefe\label{def2}
The family $\mathcal{LU}_H$ of locally univalent sense-preserving harmonic functions $f$ in the disk $\ID$ of the form \eqref{eq1}
is called a \textit{linear invariant family} {\rm (LIF)} if for each $f=h+\overline{g}\in \mathcal{LU}_H$ the following conditions
are fulfilled: $a_1=1$ and
$$\frac{f(\phi(z))-f(\phi(0))}{h'(\phi(0))\phi'(0)}\in \mathcal{LU}_H
%\;\; \forall\; \omega(z)=e^{i\theta}\frac{z+a}{1+\overline{a} z},\;a\in \ID,\;\theta\in\mathbb{R}, \omega\in {\rm Aut\,}(\ID)
$$
for each $\phi\in {\rm Aut\,}(\ID)$.
% $\omega (z)=e^{i\theta}\frac{z+a}{1+\overline{a} z}$, $a\in \ID$ and $\theta\in\mathbb{R}$.
A family $\mathcal{AL}_H$ is called \textit{linear and affine invariant} {\rm (ALIF)} if it is {\rm LIF} and in addition
each $f\in \mathcal{AL}_H$ satisfies the condition that
$$\frac{f(z)+\varepsilon \overline{f(z)}}{1+\varepsilon \overline{f_{\overline z}(0)}}\in \mathcal{AL}_H ~\mbox{ for every }~ \varepsilon\in\ID.
$$
The number $\texttt{ord}\,\mathcal{AL}_H=\sup_{f\in\mathcal{AL}_H}|a_2|$ is known as the order of the {\rm ALIF} $\mathcal{AL}_H$.
\edefe

The order of LIF $\mathcal{LU}_H$ without the assumption of affine invariance property is defined in the same way:
$\texttt{ord}\,\mathcal{LU}_H=\sup_{f\in\mathcal{LU}_H}|a_2|$.

Throughout the discussion, we suppose that the orders of these families, namely, $\texttt{ord}\,\mathcal{AL}_H$ and
$\texttt{ord}\,\mathcal{LU}_H$, are finite. The universal linear and affine invariant family, denoted by $\mathcal{AL}_H(\alpha)$,
is the largest ALIF $\mathcal{AL}_H$ of order $\alpha=\texttt{ord}\,\mathcal{AL}_H$. Thus, the subfamily $\mathcal{AL}_H^0$
of ALIF $\mathcal{AL}_H$ consists of all functions $f\in \mathcal{AL}_H$ such that $f_{\overline z}(0)=0.$
If $f\in \mathcal{AL}_H^0$ is univalent in $\ID$, then according to the result of Sheil-Small \cite{Sheil90}
one has the following sharp lower estimate:
\be\label{eq2}
d_f(0)\ge\frac{1}{2\alpha}.
%\eqno (2)
\ee

For $\alpha >0$ and $Q\geq 1$, denote by ${\mathcal H}(\alpha,Q)$ the set of all locally univalent $Q$-quasiconformal
harmonic mappings $f=h+\overline{g}$ in $\ID$ of the form
\eqref{eq1} with the normalization $a_1+a_{-1}=1$ such that
$$h(z)/h'(0)\in{\mathcal U}_\alpha,\;\;|g'(z)/h'(z)|\le k,\;k=(Q-1)/(Q+1)\in [0,1).
$$

The family ${\mathcal H}(\alpha,Q)$ was introduced and investigated in details \cite{Star93,Star95}.
In particular, he established double-sided estimates of the value $d_f(z)$ for functions belonging to the family ${\mathcal H}(\alpha,Q)$
(see \cite{S2}).

%\bdefe (cf. \cite{Star93,Star95}) \label{def3}
%Denote by ${\mathcal H}(\alpha,Q)$ the set of all locally univalent $Q$-quasiconformal harmonic mappings in the disk $\ID$ of the form
%\eqref{eq1} with the normalization $a_1+a_{-1}=1$ such that
%$$h(z)/h'(0)\in{\mathcal U}_\alpha,\;\;|g'(z)/h'(z)|\le k,\;k=(Q-1)/(Q+1)\in [0,1).
%$$
%\edefe

Note that the classes ${\mathcal H}(\alpha,Q)$, which expand with the increasing values of $\alpha\in [1,\infty]$ and $Q\in[1,\infty]$,
cover all sense-preserving locally quasiconformal harmonic mappings with the indicated normalization.

We shall restrict ourselves to the case of finite $Q$.
In \cite{Star93,Star95}, it was also shown that the family ${\mathcal H}(\alpha,Q)$ possess the property of linear invariance in the following sense:
for each $f=h+\overline{g}\in {\mathcal H}(\alpha,Q)$ and for every $\phi (z) =e^{i\theta}\frac{z+a}{1+\overline{a} z}\in {\rm Aut\,}(\ID)$,
the transformation
\be\label{eq3}
\frac{f(\phi(z))-f(\phi(0))}{\partial_\theta f(\phi(0))|\phi'(0)|}\in {\mathcal H}(\alpha,Q),
%\;\; \forall\; \omega(z)=e^{i\theta}\frac{z+a}{1+\overline{a} z},\;a\in \ID,\;\theta\in\mathbb{R},
%\eqno (3)
\ee
%belongs to ${\mathcal H}(\alpha,Q)$,
where $\partial_\theta f(z)=h'(z)e^{i\theta}+\overline{g'(z)e^{i\theta}}$ denotes the directional derivative of the complex-valued function
$f$ in the direction of the unit vector $e^{i\theta}$.
\medskip

%{\bf 1.}
In \cite{S2}, Starkov proved that for each $f\in {\mathcal H}(\alpha,Q)$ and $z\in\ID$,
$$\frac{1-|z|^2}{2\alpha Q}\max_\theta|\partial_\theta f(z)| \le  d_f(z) \le  Q (1-|z|^2)\min_\theta|\partial_\theta f(z)|
$$
which is equivalent to
\beq\label{eq4}
\frac{1-|z|^2}{2\alpha Q}\left(|h'(z)|+|g'(z)|\right)
\le d_f(z) \leq Q (1-|z|^2)\left(|h'(z)|-|g'(z)|\right),
\eeq
%\beq\label{eq4}
%\frac{1-|z|^2}{2\alpha Q}\left(|h'(z)|+|g'(z)|\right)
%& =&\frac{1-|z|^2}{2\alpha Q}\max_\theta|\partial_\theta f(z)| \nonumber\\
%&\le & d_f(z) \\
%&\le & Q (1-|z|^2)\min_\theta|\partial_\theta f(z)|\nonumber\\
%&=&Q (1-|z|^2)\left(|h'(z)|-|g'(z)|\right), \nonumber
%%\eqno (4)
%\eeq
and the lower estimate is sharp in contrast to the upper one.

One of the main aims of this article is to establish sharp estimations of the ratio $d_f(z)/d_h(z)$ for $Q$-quasiconformal
harmonic mappings $f= h+\overline{g}$. In particular, sharp upper estimate in \eqref{eq4} is obtained.
The ratio $d_f(z)/d_h(z)$ demonstrates how the radius of the largest univalent disk with the center at $h(z)$ on the
manifold $h(\ID)$ varies if we add, to the analytic function $h$,  the function $\overline{g}$.

We now state our first result.

\bthm\label{th1}
Let $f= h+\overline{g}\in {\mathcal H}(\alpha,Q)$ for some $Q\in[1,\infty]$, and $\omega(z)=g'(z)/h'(z)$ be the complex dilatation
of the mapping $f$. Then for $z\in \ID$,
\be\label{eq5}
1-k\le m\left (\frac{|\omega(z)|}{k},Q\right) \le \frac{d_f(z)}{d_h(z)} \le M\left (\frac{|\omega(z)|}{k},k\right )\le 1+k,
%\eqno (5)
\ee
where $k=(Q-1)/(Q+1)\in [0,1]$. Here the functions $M(.,k)$ and $m(.,Q)$ are defined as follows:
\beq\label{th-neweq1}
M(x,k)= \left \{\begin{array}{cl}
\ds 1+\frac{k}{x}\left\{1-\left(\frac{1}{x}-x \right)\log \left(1+x\right)\right\} & \mbox{ when }x\in(0,1]\\[4mm]
\ds \lim_{x\ra 0^+} M(x,k) =1+\frac{k}{2} & \mbox{ when } x=0
\end{array}\right . ,
\eeq
%for $x\in [0,1],\; k\in [0,1],\;Q=(1+k)/(1-k)\in[1,\infty]$
and
\beq\label{th-neweq2}
\frac{1}{m(x,Q)}= \left \{\begin{array}{cl}
\ds \int_0^1\frac{1+\varphi^{-1}(\varphi(t)/Q)x}{1-kx+\varphi^{-1}(\varphi(t)/Q)(x-k)}\,dt & \mbox{ when }Q<\infty\\[4mm]
\ds 0 & \mbox{ when } Q=\infty
\end{array}\right . ,
\eeq
with
$$\varphi(t)=\frac{\pi}{2}\frac{\K'(t)}{\K(t)} \quad (t\in(0,1))
$$
where $\K$ denotes the (Legendre) complete elliptic integral of the first kind given by
$$\K(t)= \int_0^{\pi/2} \frac{dx}{\sqrt{1-t^2\sin^2 x}} =\int_0^{1} \frac{dx}{\sqrt{(1-x^2)(1-t^2x^2)}}
$$
and $\K'(t)=\K(\sqrt{1-t^2})$. The argument $t$ is sometimes called the modulus of the elliptic integral $\K(t)$.

 Estimations in \eqref{eq5} are sharp for the family ${\mathcal H}(\alpha,Q)$ for $Q<\infty$ and for each $\alpha\ge 1$.
When $Q=\infty$, estimations in \eqref{eq5} are sharp in the sense that for each $z\in\ID$,
$$\inf_{f\in {\mathcal H}(\alpha,\infty)}\frac{d_f(z)}{d_h(z)}=m(x,\infty)=0
~\mbox{ and }~\sup_{f\in {\mathcal H}(\alpha,\infty)}\frac{d_f(z)}{d_h(z)}=M(1,1)=2.
$$
\ethm

\br\label{rem1}
For fixed $\zeta\in\ID$, the least value of the upper estimation in \eqref{eq5} is attained when $x=0$; that is
when $\omega(\zeta)=0$. In this case the estimation in \eqref{eq5} takes the form
$$\frac{d_f(\zeta)}{d_h(\zeta)}\le1+\frac{k}{2}.
$$
\er

Suppose that $f =h +\overline{g}\in {\mathcal H}(\alpha,Q),\,\alpha\in [1,\infty]$, and
$f_1(z)=C\cdot f(z)=h_1(z)+\overline{g_1(z)}$, where $C$ is a complex constant. Then
the following relations hold:
$$d_{f_1}(z)=|C|\,d_f(z) ~\mbox{ and }~d_{h_1}(z)=|C|\,d_h(z), \quad z\in\ID.
$$
Moreover, after appropriate normalization, every $Q$-quasiconformal harmonic mapping in $\ID$ belongs to  the family ${\mathcal H}(\alpha,Q)$
for some $\alpha$. Therefore an equivalent formulation of Theorem \ref{th1} may now be stated.

\bthm\label{th1a}
Let $f = h +\overline{g }$ be a locally univalent $Q$-quasiconformal harmonic mapping of the disk $\ID$, $Q\in[1,\infty]$,
and $\omega(z)=g'(z)/h'(z)$. Then the inequalities \eqref{eq5} continue to hold
%$$1-k\le m(|\mu(z)|/k,Q) \le d_f(z)/d_h(z) \le M(|\mu(z)|/k,k)\le 1+k ~\mbox{ for $z\in \ID $}
%$$
and the estimations in \eqref{eq5} are sharp.
%Here $m(x,Q)$ and $M(x,k)$ are defined as in Theorem \ref{th1}.
\ethm

Next, we consider $f= h+\overline{g}\in {\mathcal H}(\alpha,Q)$ and introduce $H(z)=h(z)/h'(0)$ from ${\mathcal U}_\alpha$.
Then we have (see \cite{Star93,Star95})
$$\frac{1}{1+k}\le|h'(0)|\le\frac{1}{1-k}
$$
and thus,
$$\frac{d_H(z)}{1+k}\le d_h(z)=|h'(0)|\cdot d_H(z)\le \frac{d_H(z)}{1-k}.
$$
These inequalities and \eqref{eq5} give the following.

%\medskip
%{\bf Consequence 1.}
\bcor
Let $f = h+\overline{g}\in {\mathcal H}(\alpha,Q)$ and $h(z)=h'(0)H(z)$. Then we have
$$\frac{d_H(z)}{Q}\le d_f(z)\le Q\, d_H(z)~\mbox{ for $z\in \ID $}.
$$
\ecor

The sharpness of the last double-sided inequalities at the point $z=0$ follows from the proof of Theorem \ref{th1}.

We now state the remaining results of the article.

\bthm \label{th2}
Let $f =h+\overline{g}$ be a locally quasiconformal harmonic mapping belonging to the family $\mathcal{AL}_H$ with
$\texttt{ord}\,(\mathcal{AL}_H)=\alpha<\infty$, $\omega(z)=g'(z)/h'(z)$ and $|\omega(z)|<1$. Then
\be\label{ext1a}
d_f(z)\ge\frac{1-|\omega(z)|}{2\alpha}\left(\frac{1-|z|}{1+|z|}\right)^\alpha ~\mbox{ for $z\in \ID $}.
\ee
The estimation $d_f(0)$ is sharp for example in the universal {\rm ALIF} $\mathcal{AL}_H(\alpha)$.
\ethm

Recall that a locally univalent function $f$ is said to be convex in the disk
$\ID(z_0,r):=\{z:\,|z-z_0|<r\}$
if $f$ maps $\ID(z_0,r)$ univalently onto a convex domain.
The radius of convexity of the family $\mathcal F$ of functions defined on the disk $\ID$ is the largest number
$r_0$ such that every function $f\in{\mathcal F}$ is convex in the disk $\ID(0,r_0)$.

\bthm\label{th3} %{\bf Theorem 3.}
 If $f\in {\mathcal H}(\alpha,Q)$,
then for every $z\in\ID$, the function $f$ is convex in the disk $\ID (z,R(z))$, where
\be\label{eq17a}
R(z)=\frac12\left(R_0+R^{-1}_0-\sqrt{\left(R_0-R^{-1}_0\right)^2+4|z|^2}\right),
\ee
and
\be\label{eq17b}
R_0=\alpha+k^{-1}-\sqrt{k^{-2}-1}-\sqrt{\left(\alpha+k^{-1}-\sqrt{k^{-2}-1}\right)^2-1}.
%\eqno (17)
\ee
In particular, the radius of convexity of the family ${\mathcal H}(\alpha,Q)$ is no less than $R_0$.
\ethm

 The proofs of Theorems \ref{th1}, \ref{th2} and \ref{th3} will be presented in Section \ref{sec2}.

\section{Proofs of the Main results}\label{sec2}

\subsection{Proof of Theorem \ref{th1}}
%\bpf 1).
The proof of the theorem is divided into three parts.\\

\noindent{\bf Part 1:} Let $f=h+\overline{g}$ satisfy the assumptions of Theorem \ref{th1}.
In compliance with the definition of the value $d_f(0)$, there exists a boundary point $A$ of the manifold $f(\ID)$ such that
$A\in\{w:\,|w|=d_f(0)\}$. Consider the smooth curve $\ell _0=f^{-1}([0,A))$, namely, the preimage of the semi-open segment
$[0,A)$ with the starting point $0$ in the disk $\ID$. Then
$$d_f(0)=|A|=\left|\int_{\ell_0}df(z)\right|=\min_{\gamma} \left|\int_{\gamma}df(z)\right|,
$$
where the minimum is taken over all smooth paths $\gamma(t),\,t\in [0,1)$, such that $\gamma(0)=0,\,|\gamma(t)|<1$
and $\lim_{t\to 1^{-}}|\gamma(t)|=1.$

Similarly we define the value
$$d_h(0)=|B|=\left|\int_{\ell}dh(z)\right|=\min_{\gamma} \left|\int_{\gamma}dh(z)\right|,
$$
where the simple smooth curve $\ell =h^{-1}([0,B))$ is emerging from the origin, the preimage of the semi-open segment
$[0,B)$ under the mapping $h$. Consider the following pa\-ra\-met\-ri\-za\-tion of the curve $\ell$:
$\ell(t)=h^{-1}(Bt),\,t\in[0,1)$. Then $h'(\ell (t))\ell'(t)=B$ and
\beq \label{eq6}
d_f(0) = \left|\int_0^1 df(\ell_0(t))\right|&\le&\left|\int_0^1df(\ell(t))\right|\nonumber\\
&=&\left|\int_0^1\!\left\{h'(\ell(t))\ell'(t)+\overline{g'(\ell(t))\ell'(t)}\right\}dt\right|\nonumber\\
&=&|B|\left|\int_0^1\!\left\{1+\frac{\overline{g'(\ell(t))\ell'(t)}}{h'(\ell(t))\ell'(t)}\right\}dt\right|\nonumber\\
&\le & d_h(0)\left\{1+\int_0^1|\omega(\ell(t))|\,dt\right\}.
%\eqno (6) Not used!
\eeq

At first we consider the case $k=\sup_{z\in\ID}|\omega(z)|<1$. Since $|\omega(z)| \le k$ for $z\in \ID$, we have
$$\omega(0)/k=\overline{a_{-1}}/(k\,a_1)=:u\in\overline{\ID}.
$$
If $|u|=1$ for $k<1$, then we have the inequality
$$d_f(0)\le d_h(0)(1+k)=d_h(0)M(1,k)
$$
which proves the upper estimate in the inequality \eqref{eq5} for $z=0$.

Let us now assume that $|u|<1$ for some $k<1$. Then, from a generalized version of the classical Schwarz lemma
(see for example \cite[Chapter VIII, \S 1]{Golu66}), it follows that
\be\label{eq7}
\frac{|\omega(z)|}{k}\le \frac{|z|+|u|}{1+|u|\,|z|}.
%\eqno (7)
\ee
Consequently, by \eqref{eq6}, one has
\be\label{eq8}
d_f(0)\le d_h(0)\left\{1+k\int_0^1\frac{|\ell(t)|+|u|}{1+|u|\, |\ell(t)|}\,dt\right\}.
%\eqno (8)
\ee

Also, the function $h^{-1}(B\zeta)$ maps biholomorphically $\ID$ onto some subdomain of the disk $\ID$.
Applying the classical Schwarz lemma, we obtain the inequality $|h^{-1}(B\zeta)|\le|\zeta|$ and hence, $|\ell(t)|\le t$ holds.
Using the last estimate and the inequality \eqref{eq8}, one can obtain, after evaluating the integral, the inequality
$$d_f(0)\le d_h(0)\left\{1+k\int_0^1\frac{t+|u|}{1+|u|t}\,dt\right\}=d_h(0)M(|u|,k),
$$
where $M(x,k)$ is defined by \eqref{th-neweq1}. The function $M(x,k)$ is strictly increasing on $(0,1]$ with respect to the variable $x$
and for each fixed $k\in[0,1]$. This follows from the observation that (see \eqref{th-neweq1})
\beqq
\frac{\partial M(x,k)}{\partial x}
%& =& \frac{\partial}{\partial x}\left\{1+\frac{k}{x}\left (1-\left(\frac{1}{x}-x \right)\log \left(1+x\right)\right )\right\}\\
&=& -\frac{k}{x^2}+\frac{2k}{x^3}\log(1+x)-k\left (\frac{1-x}{x^2}\right ),
\eeqq
which is positive, since $\log(1+x)%=x-\frac{x^2}{2}+\frac{x^3}{3}-\frac{x^4}{4}+\dots
>x-x^2/2$. Hence
\be\label{eq9}
d_f(0)\le d_h(0)M(|u|,k)\le d_h(0)M(1,k)=(1+k)d_h(0).
%\eqno (9)
\ee

We now set $k=1$. According to Lewy's theorem \cite{Le}
for locally univalent harmonic mapping $f$, we obtain that $|\omega(z)|\ne 1$ for all $z\in\ID$. Next we obtain the inequality \eqref{eq9}
in the case $k=1$ by repeating the argument of the case $k<1$.

We now begin to prove that the upper estimate in \eqref{eq5} is true for all $\zeta\in\ID.$ As mentioned above,
the family ${\mathcal H}(\alpha,Q)$ is linear invariant in the sense of \cite{Star93,Star95} (see \eqref{eq3} above).
Hence, for each fixed $\zeta=re^{i\theta}\in\ID$ $(r\in[0,1),\,\theta \in \mathbb{R})$, the function $F$ defined by
$$F(z)=\frac{f\left(e^{i\theta}\frac{z+r}{1+r z}\right)-f(re^{i\theta})}{\partial_\theta f(re^{i\theta})(1-r^2)}=H(z)+\overline{G(z)}
$$
belongs to the family ${\mathcal H}(\alpha,Q)$, where $H$ and $G$ are analytic in $\ID$ such that $H(0)=G(0)=0.$
Therefore, in view of \eqref{eq9} for $k\in[0,1]$, we have
$$d_F(0)=\frac{d_f(\zeta)}{|\partial_\theta f(\zeta)|(1-|\zeta|^2)}\le d_H(0)M(x,k),
$$
where
$x=|G'(0)/H'(0)|/k=|\omega(\zeta)|/k\in [0,1]
$
if $k\in[0,1)$, and $x=|G'(0)/H'(0)|=|\omega(\zeta)|\in [0,1)$ when $k=1$.
Note that
$$H(z)=\frac{h\left(e^{i\theta}\frac{z+r}{1+r z}\right)-h(re^{i\theta})}{\partial_\theta f(re^{i\theta})(1-r^2)}.
$$
Consequently,
$$d_H(0)=\frac{d_h(\zeta)}{|\partial_\theta f(\zeta)|(1-|\zeta|^2)}
$$
so that
$$d_f(\zeta)\le d_h(\zeta)M(x,k)\le (1+k)d_h(\zeta)
$$
and we complete the proof of the upper estimate in \eqref{eq5}. \\

\noindent{\bf Part 2:} We now deal with the sharpness of the upper estimate in \eqref{eq5}. Consider the case $k\in[0,1)$.
For every $\alpha\in\mathbb N$ and every $\zeta\in\ID$, we shall indicate functions from the families
${\mathcal H}(\alpha,Q)$ such that
$d_f(\zeta)/d_h(\zeta)=M(x)=1+k$, where $x=|\omega(\zeta)|/k.$ Since the families ${\mathcal H}(\alpha,Q)$ are enlarging
with increasing values of $\alpha$, the sharpness of the upper estimate in \eqref{eq5} will be shown for every
$\zeta\in\ID$ and each $\alpha\in[1,\infty]$.

Consider the sequence  $\{k_n\}_{n=1}^\infty$ of functions from ${\mathcal U}_n$ defined by
$$k_n(z)=\frac{i}{2n}\left[\left(\frac{1-iz}{1+iz}\right)^n-1\right].
$$
Then we have $d_{k_n}(0)=1/2n$ (see \cite{Pom64}) and observe that $k_n$ maps the unit disk $\ID$ univalently onto
the Riemann surface $k_n(\ID)$ whose boundary described by
$$\partial k_n(\ID)=\left\{\frac{i}{2n}[(i\lambda)^n-1]:\,\lambda\in\mathbb R\right\}=\left\{\frac{i}{2n}[s\,e^{\pm i\pi n/2}-1]:\,s\ge 0\right\}\
$$
consists of two rays. Then the univalent image of the disk $\ID$ under the mapping
$$f_n(z)=h_n(z)+\overline{g_n(z)}=\frac{1}{1-k}[k_n(z)-k\,\overline{k_n(z)}]\in {\mathcal H}(n,Q),\;\;k\in [0,1),
$$
($a_1=\frac{1}{1-k},\;a_{-1}=-\frac{k}{1-k}$) represents the manifold with the boundary
$$\partial f_n(\ID)=\left\{\frac{i}{2n(1-k)}[s(e^{\pm i\pi n/2}+k\,e^{\mp i\pi n/2})-1-k]:\,s\ge 0\right\},
$$
which consists of two rays parallel to the coordinate axes and arising from the point $-\frac{i}{2n}Q$.
Note that the function $f_n$ maps the semi-open segment $[0,-i)$ bijectively onto $[0,-\frac{i}{2n}Q)$ and thus, we conclude that
$$d_{f_n}(0)=\frac{Q}{2n}.
$$
This gives
$$d_{f_n}(0)=d_{k_n}(0)Q=d_{h_n}(0)(1+k),
$$
where $h_n(z)=k_n(z)/(1-k)$. The sharpness of the upper estimate in \eqref{eq5} is proved for $\zeta=0$ and $k<1$.

Next we let $0\ne\zeta\in\ID$, $k<1$, and consider a conformal automorphism $\phi(z)=(z+\zeta)/(1+\overline{\zeta} z)$ of the unit disk $\ID$.
Then the inverse mapping is given by $\phi^{-1}(z)=(z-\zeta)/(1-\overline{\zeta} z)$.
From the condition \eqref{eq3} of the linear invariance property of the family ${\mathcal H}(\alpha,Q)$, it follows that the function
$f$ defined by
$$f(z)=\frac{f_n(\phi^{-1}(z))-f_n(-\zeta)}{\partial_0 f_n(-\zeta)(1-|\zeta|^2)}=h(z)+\overline{g(z)}
$$
belongs to ${\mathcal H}(\alpha,Q)$, where $h$ and $g$ have the same meaning as above.
Taking into account of the normalization condition for functions in the family ${\mathcal H}(\alpha,Q)$,
we deduce that
$$\frac{f(\phi(z))-f(\zeta)}{\partial_0 f(\zeta)(1-|\zeta|^2)}=f_n(z)=h_n(z)+\overline{g_n(z)}.
$$
Therefore,
$$d_{f_n}(0)=\frac{d_f(\zeta)}{|\partial f_0(\zeta)|(1-|\zeta|^2)}=d_{h_n}(0)(1+k).
$$
On the other hand, a direct computation gives
$$h_n(z)=\frac{h(\phi(z))-h(\zeta)}{\partial_0 f(\zeta)(1-|\zeta|^2)}
~\mbox{ and }~d_{h_n}(0)=\frac{d_h(\zeta)}{|\partial f_0(\zeta)|(1-|\zeta|^2)}
$$
showing that
$$d_{f_n}(0)|\partial f_0(\zeta)|(1-|\zeta|^2)=d_f(\zeta)=d_{h}(\zeta)(1+k),
$$
which completes the proof of the upper estimation in Theorem \ref{th1} for $k\in[0,1)$.

If $k=1$ then for $j\in\mathbb N$, we consider the sequence $\{f_{n,j}\}$ of functions
$$f_{n,j}(z)=h_{n,j}(z)+\overline{g_{n,j}(z)}=j k_n(z)-(j-1)\,\overline{k_n(z)}.
$$
We see that $f_{n,j}\in {\mathcal H}(n,2j-1)\subset {\mathcal H}(n,\infty)$ for each $j\in \mathbb N$. Therefore,
$$d_{f_{n,j}}(0)=d_{h_{n,j}}(0)M(1,1-1/j).
$$
Hence
$$\sup_{j\in\mathbb N} \frac{d_{f_{n,j}}(0)}{d_{h_{n,j}}(0)}=M(1,1)=2.
$$
The sharpness of the upper estimation in \eqref{eq5} for $k=1,$ $\zeta\ne 0$, can be proved analogously. So, we omit the details.

%\begin{figure}
%\begin{center}
%\includegraphics[width=6.5cm]{fig1} %{Covering_h2.eps}
%\hspace{.3cm}
%\includegraphics[width=6.5cm]{fig2} %{Covering_f2.eps}\\
%\end{center}
%\caption{The covering disks for functions $h_2$ and $f_2$ for $k=0.25$ with centres at origin and radii $1/3$ and $5/12$ correspondingly \label{fig1-2}}
%\end{figure}
\begin{figure}
\begin{center}
\includegraphics[width=6.8cm]{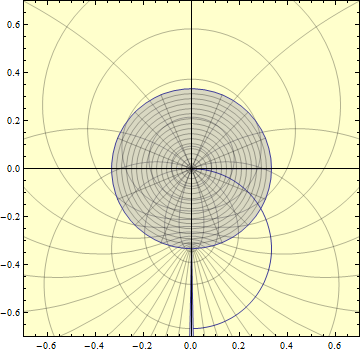} %{fig1.eps}%{fig1.png} %{Covering_h2.eps}
\end{center}
\textbf{(a)} $h_2$
\hspace{.3cm}
\begin{center}
\includegraphics[width=6.8cm]{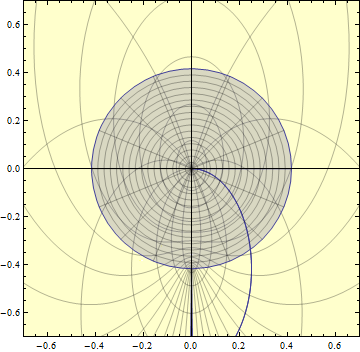} %{fig2.eps}%{fig2.png} %{Covering_f2.eps}\\
\end{center}
\textbf{(b)} $f_2$
\caption{The covering disks for functions $h_2$ and $f_2$ for $k=0.25$ with centers at origin and radii $1/3$ and $5/12$, respectively.  \label{fig1-2}}
\end{figure}
The images of polar grid in the unit disk under mappings $h_2$ and $f_2$ are indicated in Figures \ref{fig1-2}(a)-(b) which
illustrate the sharpness assertions proved in the above estimations.\\

\noindent{\bf Part 3:} Finally, we deal with lower estimation of $d_f(z)$. If $k=1$, then the lower
estimation in \eqref{eq5} is trivial because $m(x,\infty)=0$. So, we may assume that $k\in[0,1)$.
As in Part~1, we define the boundary points $A$ and $B$ of the manifolds $f(\ID)$ and $h(\ID)$, respectively,
and smooth curves $\ell_0=f^{-1}([0,A))$ and $\ell$ in the same manner as in Part~1.
Consider the parametrization of the curve $\ell_0$:
$$\ell_0(t)=f^{-1}(At),\,t\in[0,1).
$$
Then $df(\ell_0(t))=A dt$ and thus,
\beq\label{eq10}
d_h(0)&=&\left|\int_0^1 dh(\ell(t))\right|\nonumber \\
&\le &\left|\int_0^1 h'(\ell_0(t))\ell_0'(t)\,dt\right|\nonumber\\
&=&\left|\int_0^1\!\left\{h'(\ell_0(t))\ell_0'(t)+\overline{g'(\ell_0(t))\ell_0'(t)}\right\}\right. \nonumber \\
&& \times \left .\left(1-\frac{\overline{g'(\ell_0(t))\ell_0'(t)}}{h'(\ell_0(t))\ell_0'(t)+\overline{g'(\ell_0(t))\ell_0'(t)}}\right)dt\right| \nonumber\\
&=& \left|\int_0^1\frac{h'(\ell_0(t))\ell_0'(t)}{h'(\ell_0(t))\ell_0'(t)+\overline{g'(\ell_0(t))\ell_0'(t)}}df(\ell_0(t))\right| \nonumber\\
&\le & |A|\int_0^1\frac{dt}{1-|\omega(\ell_0(t))|}.
\eeq
%\eqno (10)

In view of the inequality \eqref{eq7}, we find that
\be\label{eq11}
|\omega(\ell_0(t))|\le k\frac{|\ell_0(t)|+x}{1+x|\ell_0(t)|},
%\eqno (11)
\ee
where $x=|\omega(0)|/k$.

It is possible to obtain an estimate for $|\ell_0(t)|=|f^{-1}(At)|,\,t\in [0,1)$, with the help of the analog of the Schwarz lemma for
$Q$-quasiconformal automorphisms of the disk. Let $F$ be a $Q$-quasiconformal automorphism of $\ID$, and $F(0)=0$. It is known
(see for example \cite[Chapter 10, equality (10.1)]{AVV97}) that the sharp estimation
$$|F(z)|\le \varphi^{-1}\left(Q^{-1}\varphi(|z|)\right)
$$
holds, where $\varphi$ and $Q$ are as in the statement.
%$\varphi(t)=\frac{\pi}{2}\frac{K'(t)}{K(t)},\;\;K(t),\,K'(t)$ -- complete elliptic integrals of Jacobi, $t\in(0,1)$.
The function $f^{-1}(A w)$ defined on the unit disk $\{w:\,|w|<1\}$ satisfies the conditions $f^{-1}(0)=0$ and $|f^{-1}(A w)|<1$. Let $\Phi$ be
the univalent conformal mapping of the domain $f^{-1}(A\ID)$ onto the unit disk $\ID$ and $\Phi(0)=0.$
Then the composition $\Phi\circ f^{-1}(Az)$ is a $Q$-quasiconformal automorphism of $\ID$ and $\Phi^{-1}$ satisfies the conditions
of the classical Schwarz lemma for analytic functions. Hence, we have
$$|\ell_0(t)|=|\Phi^{-1}(\Phi\circ f^{-1}(At))|\le|\Phi\circ f^{-1}(At)|\le \varphi^{-1}(Q^{-1}\varphi(t)).
$$
As a result of it and taking into account of the last estimation, inequalities \eqref{eq10} and \eqref{eq11}, and the fact that the function
$(1+y\,x)/(1-kx+y(x-k))$ is strictly increasing with respect to $y$ on $(0,1)$, we conclude that
$$d_h(0)\le d_f(0)\int_0^1\frac{1+y\,x}{1-kx+y(x-k)}\,dt\le \frac{d_f(0)}{1-k},
$$
where $y=\varphi^{-1}(Q^{-1}\varphi(t))\le 1$ for $t\in(0,1)$.
Therefore the lower estimate in \eqref{eq5} is sharp at the origin.

The proof of the lower estimation in \eqref{eq5} for $0\ne\zeta\in \ID$ follows easily if we proceed with the
same manner as in Part~1 and use the linear invariance property of the family ${\mathcal H}(\alpha,Q)$.

For the sharpness of the left side of the inequality in \eqref{eq5} for $k\in[0,1)$, we consider the functions (see \cite{Star93,Star95})
\be\label{eq12}
h_\alpha(z)=\frac{1}{2i\alpha}\left[\left(\frac{1+iz}{1-iz}\right)^\alpha-1\right]\in{\mathcal U}_\alpha
%\eqno (12)
\ee
and
$$f(z)=h(z)+\overline{g(z)}:= \frac{h_\alpha(z)}{1+k}+\frac{k\overline{h_\alpha(z)}}{1+k}.
$$
Then it is a simple exercise to see that
$$d_f(0)=\frac{1}{2\alpha Q} ~\mbox{ and }~ d_h(0)=\frac{1}{2\alpha(1+k)}.
$$
%\begin{figure}
%\begin{center}
%\includegraphics[width=6.5cm]{fig3} \hspace{.3cm} %{Covering_h2.eps}
%\includegraphics[width=6.5cm]{fig4} %{Covering_f2.eps}\\
%\end{center}
%\caption{The covering disks for functions $h_\alpha/(1+k)$ and $f=(h_\alpha+k\overline h_\alpha)/(1+k)$ for $k=0.25$ with centres at origin and radii $1/5$ and $3/20$ correspondingly \label{fig3-4}}
%\end{figure}
%\begin{figure}
%\begin{center}
%\includegraphics[width=6.8cm]{fig3.pdf}%{fig3.eps}%{fig3.png}
%\hspace{.3cm} %{Covering_h2.eps}
%%\textbf{(a)} $h_\alpha/(1+k)$
%\end{center}
%\textbf{(a)} $h_\alpha/(1+k)$
%\hspace{.3cm}
%\begin{center}
%\includegraphics[width=6.8cm]{fig4.pdf} %{fig4.eps}%{fig4.png} %{Covering_f2.eps}\\
%\end{center}
%\textbf{(b)} $f=(h_\alpha+k\overline h_\alpha)/(1+k)$
%\caption{Covering disks for functions $h_\alpha/(1+k)$ and $f=(h_\alpha+k\overline h_\alpha)/(1+k)$ for $k=0.25$
%with centers at origin and radii $1/5$ and $3/20$, respectively.    \label{fig3-4}}
%\end{figure}
\begin{figure}
\begin{center}
\includegraphics[width=6.8cm]{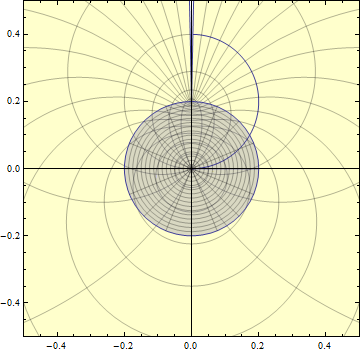}%{fig3.eps}%{fig3.png}
\hspace{.3cm} %{Covering_h2.eps}
%\textbf{(a)} $h_\alpha/(1+k)$
\end{center}
\textbf{(a)} $h_\alpha/(1+k)$
\hspace{.3cm}
\begin{center}
\includegraphics[width=6.8cm]{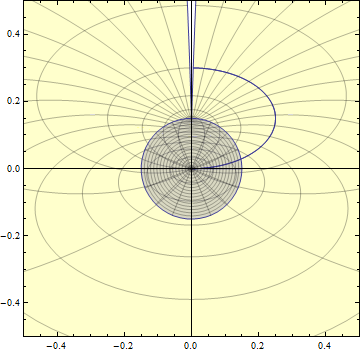} %{fig4.eps}%{fig4.png} %{Covering_f2.eps}\\
\end{center}
\textbf{(b)} $f=(h_\alpha+k\overline h_\alpha)/(1+k)$
\caption{Covering disks for functions $h_\alpha/(1+k)$ and $f=(h_\alpha+k\overline h_\alpha)/(1+k)$ for $k=0.25$
with centers at origin and radii $1/5$ and $3/20$, respectively.    \label{fig3-4}}
\end{figure}
Comparison of radii $d_h(0)$ and $d_f(0)$ and sharpness of the lower estimation of $d_f(0)/d_h(0)$ is illustrated
in Figures \ref{fig3-4}(a)--(b). In these figures, the images of polar grid in the unit disk under mappings $h_\alpha/(1+k)$ and $f$ are indicated.

If $k\to 1^{-}$ then from the last equality we obtain
$$\lim_{k\to 1^{-}} d_f(0)=0~\mbox{ and }~\lim_{k\to 1^{-}} d_h(0)=\frac{1}{4\alpha},
$$
so that
$$\inf_{f\in {\mathcal H}(\alpha,\infty)} \frac{d_f(0)}{d_h(0)}=0.
$$
Thus the last equality is sharp not only at the origin but also at points $z\in\ID$, in view of the degeneration of functions
$$f(z)=\frac{h_\alpha(z)+k\overline{h_\alpha(z)}}{1+k}
$$
when $k\to 1^{-}$. The proof of the theorem is complete.
\hfill $\Box$

\br\label{rem2}
In some neighbourhood of the origin, it is also possible to obtain a simple lower estimate in the inequality \eqref{eq5}
without the involvement of elliptic integrals. For example, the well-known theorem of Mori \cite{Mori56}
reveals that for $Q$-quasiconformal automorphism $F$ of the disk $\ID$ such that $F(0)=0$, one has
$$|F(z)|\le 16|z|^{1/Q}.
$$
Using this result in the estimation of $|\ell_0(t)|$ in Part~3 of the proof of Theorem \ref{th1}, one can easily obtain that
$$|\ell_0(t)|=|f^{-1}(At)|=\left\{
\aligned
16 t^{1/Q} & \text{ when } 0\le t< 1/16^{-Q}\\
1 & \text{ when } 16^{-Q}\le t< 1.
\endaligned
\right.
$$
The last relation provides an opportunity to estimate the ratio $d_h(z)/d_f(z)$ by means of an
integral of an elementary function, namely,
\beqq
\frac{d_h(z)}{d_f(z)}&\le &\frac{1}{m(x,Q)}\\
&\le & \frac{1}{1-k} \left(1-16^{-Q}+(1-k)\int_0^{16^{-Q}}
\frac{1+y\,x}{1-kx+y(x-k)}\,dt\right)\\
&\le &\frac{1}{1-k},
\eeqq
where $y=16t^{1/Q}\le 1$ for $t\in[0,16^{-Q}].$ Here $x=|\omega(z)|/k,\;z\in\ID$.
\er

\subsection{Proof of Theorem \ref{th2}} %\bpf
We first prove the inequality \eqref{ext1a} for $z=0$. As in the proof of Theorem \ref{th1},
consider on the circle $\{w:\,|w|=d_f(0)\}$ the boundary point
$A$ of the manifold $f(\ID)$ and define a curve $\ell_0=f^{-1}([0,A))$ with the starting point $0$ in $\ID$. Then
\be\label{eq13}
d_f(0)=|A|=\left|\int_{\ell_0}df(\zeta)\right|=\int_{\ell_0}|df(\zeta)|\ge\int_{\ell_0}\left(|h'(\zeta)|-|g'(\zeta)|\right)|d\zeta|.
%\eqno (13)
\ee

In view of the affine invariance property of the family $\mathcal{AL}_H$, the
function $F$ defined by
$$ F(\zeta)=H(\zeta)+\overline{G(\zeta)}=\frac{f(\zeta)-\varepsilon \overline{f(\zeta)}}{1-\varepsilon \overline{f_{\overline z}(0)}}
$$
belongs $\mathcal{AL}_H$ for every $\varepsilon$ with $|\varepsilon|<1$ .

For a fixed $\zeta$, we introduce $\theta(z)=\arg h'(z)-\arg g'(z)$ when $g'(z)\ne 0$,
and $\theta(z)=\arg h'(z)$ otherwise. Consider then $\varepsilon=s e^{i\theta(z)}$ for $s\in[0,1)$.
Therefore, taking into account of the relation $\overline{f_{\overline z}(0)}=\omega(0)$, we obtain that
$$H'(\zeta)=\frac{h'(\zeta)-s g'(\zeta)e^{i\theta(z)}}{1-\varepsilon \omega(0)}
$$
and thus,
\be\label{eq14}
|H'(\zeta)|\le \frac{|h'(\zeta)|-s |g'(\zeta)|}{1-s |\omega(0)|}.
%\eqno (14)
\ee
For the other side of the inequality for functions in the family $\mathcal{AL}_H$, the inequality
\be\label{eq15}
|H'(\zeta)|\ge\frac{(1-|\zeta|)^{\alpha-1}}{(1+|\zeta|)^{\alpha+1}}
%\eqno (15)
\ee
holds, where $\alpha=\texttt{ord}\,(\mathcal{AL}_H)$ is defined as in the sense of Definition \ref{def2}.
The inequality \eqref{eq15} was obtained in \cite{Sheil90} for ALIF of univalent harmonic mappings,
but the proof is still valid without a change for any ALIF $\mathcal{AL}_H$ of finite order $\alpha$.
Using inequalities \eqref{eq14} and \eqref{eq15}, we obtain the inequality
$$|h'(\zeta)|-s|g'(\zeta)|\ge (1-s |\omega(0)|)\frac{(1-|\zeta|)^{\alpha-1}}{(1+|\zeta|)^{\alpha+1}}
$$
for every $s\in(0,1)$. Allowing in the last inequality $s\to 1^{-}$ and substituting the resulting estimate
into \eqref{eq13}, we easily obtain that
\beqq
d_f(0)&\ge& (1-|\omega(0)|)\int_{\ell_0}\frac{(1-|\zeta|)^{\alpha-1}}{(1+|\zeta|)^{\alpha+1}}|d\zeta| \\
&\ge & (1-|\omega(0)|)\int_0^1\frac{(1-t)^{\alpha-1}}{(1+t)^{\alpha+1}}\,dt %\\ &=&
= \frac{1-|\omega(0)|}{2\alpha}.
\eeqq

If $0<|z|<1,$ then as in the proof of Theorem \ref{th1}, we may use the linear invariance
property of the family $\mathcal{AL}_H$ in accordance with the function $F_1\in \mathcal{AL}_H$, where
$$F_1(\zeta)=\frac{f\left(\frac{\zeta+z}{1+\overline{z} \zeta}\right)-f(z)}{h'(z)(1-|z|^2)}.
$$
In this way, applying the estimation of $d_f(0)$ to the function $F_1$, we see that
$$d_{F_1}(0)\ge\frac{1-|\omega_{f}(0)|}{2\alpha}.
$$
Also, we have
$$d_{F_1}(0)=\frac{d_f(z)}{|h'(z)|(1-|z|^2)}.
$$
It remains to note that $|\omega_{F_1}(0)|=|\omega_f(z)|$ and apply the inequality \eqref{eq15} to the function $h'(z)$.

In order to prove the sharpness of the estimate of $d_f(0)$, we first note that the functions
$p(z)=h_\alpha(z)+k\overline {h_\alpha(z)}$, where each $h_\alpha$ has the form \eqref{eq12},
belong to $\mathcal{AL}_H(\alpha)$ for every $k=|\omega(0)|\in [0,1)$. Indeed, for each $\alpha$, the
function $p$ is locally univalent and meet the normalization condition of the family
$\mathcal{AL}_H(\alpha)$, and $|p_{zz}(0)/2|=|h_\alpha''(0)/2|=\alpha$.
Affiliation of the functions
$$q(z)=\frac{p(\phi(z))-p(\phi(0))}{h_\alpha'(\phi(0))\phi'(0)}=\frac{h_\alpha(\phi(z))-h_\alpha(\phi(0))}{h_\alpha'(\phi(0))\phi'(0)}
+k\frac{\overline{h_\alpha(\phi(z))-h_\alpha(\phi(0))}}{h_\alpha'(\phi(0))\phi'(0)},
$$
and
\beqq
w(z)&=&\frac{q(z)+\varepsilon \overline{q(z)}}{1+\varepsilon \overline{q_{\overline z}(0)}}\\
&=& \frac{h_\alpha(\phi(z))-h_\alpha(\phi(0))}{h_\alpha'(\phi(0))\phi'(0)} \\
&& + \frac{\overline{h_\alpha(\phi(z))-h_\alpha(\phi(0))}}{h_\alpha'(\phi(0))\phi'(0)}
\left (\frac{k+\varepsilon h_\alpha'(\phi(0))\phi'(0)/\overline{(h_\alpha'(\phi(0))\phi'(0))}}{1+\varepsilon k\,h_\alpha'(\phi(0))\phi'(0)/\overline{(h_\alpha'(\phi(0))\phi'(0))}} \right )
\eeqq
to the family $\mathcal{AL}_H(\alpha)$ for every conformal automorphism
$\phi$ of the disk $\ID$ and every $\varepsilon\in\ID$, follow from the membership of the function $h_\alpha$ to the universal LIF ${\mathcal U}_\alpha$.
The analogous reasoning is true after the change of order of the linear and affine transforms of the function $p$.

Therefore, $p=h_\alpha +k\overline {h_\alpha}\in \mathcal{AL}_H(\alpha)$ for each $k\in [0,1)$ and at the same time
$$d_p(0)=\frac{1-|\omega(0)|}{2\alpha},
$$
which proves the sharpness of the established estimate in the universal ALIF $\mathcal{AL}_H(\alpha)$.
The proof of the theorem is complete.
\hfill $\Box$
%\epf

\br\label{rem3}
Recall that a domain $D\subset\mathbb{C}$ is called
close-to-convex if its complement $\mathbb C\setminus D$ can be written as an union of disjoint rays or lines.
The family ${\mathcal C}_H$ of all univalent sense-preserving harmonic mappings $f$ of the form \eqref{eq1}
such that $a_1=1$ and $f(\ID)$ is close-to-convex, is {\rm ALIF} (cf. \cite{Sheil90}). Also, the inequality in Theorem \ref{th2} is
sharp in the {\rm ALIF} ${\mathcal C}_H$. The order of the family ${\mathcal C}_H$ is proved to be $3$ (\cite{Clunie-Small-84}).
The harmonic analog of the analytic Koebe function $k(z)=z/(1-z)^2$ (see for example \cite[Chapter 5, p. 82]{Duren:Harmonic})
is given by
$$F(z)=\frac{z-\frac12z^2+\frac16z^3}{(1-z)^3}+\overline{\left(\frac{\frac12z^2+\frac16z^3}{(1-z)^3}\right)},
$$
where $F\in {\mathcal C}_H$ and $F(\ID)=\mathbb C\setminus (-\infty,-1/6]$ which is indeed a domain
starlike with respect to the origin. From the affine invariance property of the family ${\mathcal C}_H$,
we deduce that for every $b\in [0,1)$, the affine mapping
$$f(z)=F(z)-b\,\overline{F(z)}
$$
belongs to ${\mathcal C}_H$ such that $\omega(0)=f_{\overline z}(0)/f_z(0)=-b.$
The function $f$ is a composition of the univalent harmonic mapping $F$ of the disk $\ID$ onto
$\mathbb C\setminus (-\infty,-1/6]$ and affine transformation $\psi(w)= w-b\overline w$. The plane $\mathbb C$
with a slit $(-\infty,-1/6]$ under the transformation $\psi$ is the plane with a slit along the
ray emanating from the point $\psi(-1/6)=-(1-b)/6$ through the point $\psi(-1)=b-1<(b-1)/6$, since $b\in [0,1)$.
Therefore, $f(\ID)=\mathbb C\setminus (-\infty,-(1-b)/6]$ and thus, $d_f(0)=(1-b)/6$ and the lower estimate of
$d_f(0)$ is sharp in the {\rm ALIF} ${\mathcal C}_H$.
\er

\medskip

In the first part of the present paper, we concerned with the question about the covering of the manifold
$f(\ID)$ by disks. Now we turn our attention on the problem related with the covering of $f(\ID)$ by convex domains.

 Sheil-Small \cite{Sheil90} proved that the radius of convexity of the univalent subfamily of the linear and affine invariant family
$\mathcal{AL}_H$ of harmonic mappings is equal to
\be\label{eq16}
r_0=\alpha-\sqrt{\alpha^2-1},
%\eqno (16)
\ee
where $\alpha=\texttt{ord}\,(\mathcal{AL}_H)$. Later this result was generalized to the families of locally univalent
harmonic mappings \cite{GrafEye2010}. Now we will show the radius of convexity will be altered
under the assumption of $Q$-quasiconformality of functions $f$.

\blem\label{lem1}
Let $\mathcal{LU}_H(\alpha, Q)$ denote the {\rm LIF} of locally univalent $Q$-quasiconformal harmonic mappings of the order
$\alpha<\infty$, where $Q\le\infty$. Then the affine hull
$$ \mathcal{AL}_H=\left\{F(z)=\frac{f(z)+\varepsilon\overline{f(z)}}{1+\varepsilon \overline{ a_{-1}}} : \,f\in \mathcal{LU}_H(\alpha, Q),
\;\varepsilon\in \ID\right\}
$$
of the family $\mathcal{LU}_H(\alpha, Q)$ is linear and affine invariant of order no greater
than $\alpha+\frac{1-\sqrt{1-k^2}}{k}$, where $k=(Q-1)/(Q+1)$.
\elem\bpf
In \cite{SoStarSzy2011}, it was shown that the affine hull of the linear invariant in the sense of
Definition \ref{def2} of the family of the locally univalent harmonic mappings is the ALIF $\mathcal{AL}_H$.
Thus, it remains to determine the estimate of the order of the family $\mathcal{AL}_H$.

We begin with $F=H+\overline{G}\in\mathcal{AL}_H$. Then there exists an $f=h+\overline g\in \mathcal{LU}_H(\alpha, Q)$ of the form
\eqref{eq1} with the additional normalization $f_z(0)=h'(0)=a_1=1,$ and $\varepsilon \in\ID$ such that
$$F(z)=\frac{f(z)+\varepsilon\overline{f(z)}}{1+\varepsilon g'(0)}=H(z)+\overline{G(z)}.
$$
It is easy to compute that
$$A_2=\frac{H''(0)}{2}=\frac{a_2+\varepsilon {\overline a_{-2}}}{1+\varepsilon g'(0)},
$$
where $a_2=h''(0)/2$ and $a_{-2}=\overline{g''(0)}/2$. Taking into account of the relation $g'(z)=\omega(z)h'(z)$,
where $\omega$ is the complex dilatation of $f$ with $|\omega(z)|<k$, we see that
$$g'(0)=\omega(0) ~\mbox{ and }~g''(0)=h''(0)\omega(0)+h'(0)\omega'(0),
$$
so that
$$\overline {a_{-2}}=a_2\omega(0)+\omega'(0)/2.
$$
If we apply the Schwarz-Pick lemma (see for example \cite[Chapter VIII, \S 1]{Golu66})
to the function $\omega(z)/k$, then the inequality \eqref{eq2} in this case leads to
$$\frac{|\omega'(0)|}{k}\le 1-\frac{|\omega(0)|^2}{k^2}.
$$
Using the expression for $a_2$, we deduce that
\beqq
|A_2|&=& \left |\frac{a_2(1+\varepsilon \omega(0))+\varepsilon\omega'(0)/2}{1+\varepsilon \omega(0)}\right |\\
%&\le & |a_2|+k\frac{|\varepsilon|}{2}\frac{1-|\mu(0)/k|^2}{1-|\varepsilon \mu(0)|}\\
&\le&|a_2|+\frac{k}{2}\frac{1-|\omega(0)/k|^2}{1-|\omega(0)|}\\
&=&|a_2|+\frac{k^2-|\omega(0)|^2}{2k(1-|\omega(0)|)}
\eeqq
%\beqq
%|A_2|&= &\left |\frac{a_2(1+\varepsilon \mu(0))+\varepsilon\mu'(0)/2}{1+\varepsilon \mu(0)}\right |\\
%&\le & |a_2|+k\frac{|\varepsilon|}{2}\frac{1-|\mu(0)/k|^2}{1-|\varepsilon \mu(0)|}\\
%&\le& |a_2|+\frac{k}{2}\frac{1-|\mu(0)/k|^2}{1-|\mu(0)|}.
%\eeqq
(since $|\varepsilon|<1$). Calculating the maximum of the function $u(t)=(k^2-t^2)/(1-t)$ over the interval $[0,k]$,
we obtain the estimate
$$|A_2|\le|a_2|+\frac{1-\sqrt{1-k^{2}}}{k}\le \alpha+\frac{1-\sqrt{1-k^{2}}}{k} <\alpha+1.
$$
The proof of the lemma is complete.
\epf

Using Lemma \ref{lem1} and the equality \eqref{eq16}, one obtains the estimate of the radius of convexity of functions
in the family ${\mathcal H}(\alpha,Q)$.

\subsection{Proof of Theorem \ref{th3}}%\bpf
Let $f_0=h_0+\overline{g_0}\in {\mathcal H}(\alpha,Q)$. It is easy to see that the function $f_0$ is
convex in the same disks as the normalized function
$$f(z)=f_0(z)/h_0'(0) =h(z)+\overline{g(z)}
$$
that belongs to some LIF $\mathcal{LU}_H(\alpha,Q)$. So it is enough to prove the theorem for such functions $f$.
We first show that the function $f$ is convex in the disk centered at
the origin with radius $R_0$ defined by \eqref{eq17b}.

Clearly, the function $f$ belongs to the affine hull $\mathcal{AL}_H$ of the family
 $\mathcal{LU}_H(\alpha,Q)$.
In view of Lemma \ref{lem1}, the family $\mathcal{AL}_H$ has the order $\alpha_1\le\alpha+\frac{1-\sqrt{1-k^2}}{k}.$
Taking into consideration of the equality \eqref{eq16}, we conclude that the function $f$ is convex in the disk
of radius $R_0=\alpha_1-\sqrt{{\alpha_1}^2-1}$ centered at the origin.

We now let $0\neq z_0\in \ID$. Consider a conformal automorphism $\Phi$ of the unit disk $\ID$ given by
$$\Phi(\zeta)=e^{i\arg z_0}\left (\frac{\zeta+|z_0|}{1+|z_0|\zeta}\right ).
$$
We see that $\Phi$ maps the disk $\ID(0,R_0)$ onto the disk $\ID(z_0,R(z_0)),$ where $R(z_0)$ is
defined in \eqref{eq17a}. In view of the linear invariance property of the family $\mathcal{LU}_H(\alpha,Q)$, the function
$F$ defined by
$$F(\zeta)=\frac{f(\Phi(\zeta))-f(z_0)}{h'(z_0))\Phi'(0)}
$$
belongs to $\mathcal{LU}_H(\alpha,Q)$ and as remarked above, the function $F$ maps the disk $\ID(0,R_0)$ onto a convex domain.
Therefore, the function
$$f(z)=F(\Phi^{-1}(z))\cdot h'(z_0))\Phi'(0)+f(z_0)
$$
is convex and univalent in the disk $\ID(z_0,R(z_0))$. The proof of the theorem is complete.
\hfill $\Box$%\epf

%It is easy to note that after appropriate normalization every $f\in {\mathcal H}(\alpha,Q)$ belongs to the same linear
%invariant family in the sense of Definition \ref{def2} for the family $\mathcal{LU}_H(\alpha,Q)$ of locally univalent
%$Q$-quasiconformal harmonic mappings. As a consequence of this observation, we obtain the following.

%{\bf Consequence 3.}
%\bcor
%Let $f\in {\mathcal H}(\alpha,Q)$ for some $\alpha<\infty$ and $Q\le\infty$.
%Then for each $z\in\ID$, the function $f$ is convex and univalent in the disk $\ID(z,R(z))$,
%where $R(z)$ is defined by \eqref{eq17a}.
%\ecor

\subsection*{\bf Acknowledgements}
%\begin{acknowledgement}
The authors thank the referee for useful comments.
The research was supported by the project RUS/RFBR/P-163 under
Department of Science \& Technology (India) and the Russian Foundation for Basic Research (project 14-01-92692). The second author is currently on leave from Indian Institute of Technology Madras, India. The third author is also supported by Russian Foundation for Basic Research (project 14-01-00510) and
the Strategic Development Program of Petrozavodsk State University.

%\pagebreak


\begin{thebibliography}{99}

\bibitem{AVV97} G.~D.~Anderson, M.~K.~Vamanamurthy, and M.~Vuorinen,
\emph{Conformal invariants, inequalities, and quasiconformal maps},
John Wiley \& Sons, 1997.

\bibitem{Bib} L. Bieberbach,
\textrm{\"Uber die Koeffizienten derjenigen Potenzreihen, welche eine schlichte Abbildung des Einheitskreises vermitteln},
\textit{Sitzungsber. Preuss. Akad. Wiss.} \textbf{38}(1916),  940--955.

\bibitem{HG}  H. Chen, P. M. Gauthier and W. Hengartner,
\textrm{Bloch constants for planar harmonic mappings,}
\textit{Proc. Amer. Math. Soc.} {\bf 128}(2000), 3231--3240.

\bibitem{CPR-2014} Sh. Chen, S. Ponnusamy and  A. Rasila,
Coefficient estimates, Landau's theorem and Lipschitz-type spaces on planar harmonic mappings,
\textit{J. Aust. Math. Soc.} \textbf{96}(2014), 198--215.

\bibitem{CPW-2011} Sh. Chen, S. Ponnusamy and  X. Wang,
Bloch and Landau's theorems for planar $p$-harmonic mappings,
\textit{J. Math. Anal. Appl.} \textbf{373}(2011), 102--110.

\bibitem{CPW1} Sh. Chen, S. Ponnusamy and  X. Wang,
Landau's theorem and Marden constant for harmonic $\nu$-Bloch mappings,
\textit{Bull. Aust. Math. Soc.} {\bf 84}(2011),  19--32.

\bibitem{CPW1a} Sh. Chen, S. Ponnusamy and  X. Wang,
Properties of some classes of planar harmonic and planar biharmonic mappings,
\textit{Complex Anal. Oper. Theory,} \textbf{5}(2011), 901--916.

\bibitem{CPW1b} Sh. Chen, S. Ponnusamy and  X. Wang,
Coefficient estimates and Landau-Bloch's theorem for planar harmonic mappings,
\textit{Bull. Malaysian Math. Sciences Soc.} \textbf{34}(2)(2011), 255--265.

\bibitem{CPW-2013a} Sh. Chen, S. Ponnusamy and  X. Wang,
Covering and distortion theorems for planar harmonic univalent mappings,
\textit{Arch. Math. (Basel)}, \textbf{101}(2013), 285--291.

\bibitem{CPW-2013b} Sh. Chen, S. Ponnusamy and  X. Wang,
Harmonic mappings in Bergman spaces,
\textit{Monatsh. Math.} \textbf{170}(2013), 325--342.

\bibitem{Clunie-Small-84} J. G. Clunie and T. Sheil-Small,
\textrm{Harmonic univalent functions},
\textit{Ann. Acad. Sci. Fenn. Ser. A I Math.} {\bf 9}(1984), 3--25.

\bibitem{Du} P.~L.~Duren,
\emph{Univalent functions} (Grundlehren der
mathematischen Wissenschaften 259, New York, Berlin, Heidelberg, Tokyo), Springer-Verlag, 1983.

\bibitem{Duren:Harmonic} P.~Duren,
\emph{Harmonic mappings in the plane},
Cambridge Tracts in Mathematics, \textbf{156}, Cambridge Univ. Press, Cambridge, 2004.

\bibitem{Golu66} G. M. Goluzin,
\textit{Geometric theory of functions of a complex variable},
Translations of Mathematical Monographs, Vol. 26, American Mathematical Society, Providence, R.I. 1969 vi+676 pp.

\bibitem{GrafEye2010} S. Yu. Graf and O.~R.~Eyelangoli,
\textrm{Differential inequalities in linear and affine-invariant families of harmonic mappings},
\textit{Russian Math. (Iz. VUZ)} \textbf{54}(10)(2010), 60--62.

\bibitem{Ko} P. Koebe,
\textrm{\"{U}ber die Uniformisierung beliebiger analytischer Kurven II},
%\textit{Nachr. Acad. Wiss. G\"{o}ttingen} (1907), 633--669.
\textit{Nachrichten von der Gesellschaft der Wissenschaften zu G\"{o}ttingen, Math.-Phys. Kl.} (1907), 633--669.
%Nachrichten von der Gesellschaft der Wissenschaften zu G\"{o}ttingen, Mathematisch-Physikalische Klasse

\bibitem{Le} H. Lewy,
On the non-vanishing of the Jacobian in certain one-to-one mappings,
\textit{Bull. Amer. Math. Soc.} \textbf{42}(1936), 689--692.

\bibitem{Minda90} D. Minda,
\textrm{The Bloch and Marden constants},
Computational methods and function theory (Valparaíso, 1989), 131--142,
Lecture Notes in Math., 1435, Springer, Berlin, 1990.

\bibitem{Mori56} A. Mori,
\textrm{On an absolute constant in the theory of quasi-conformal mappings},
\textit{J. Math. Soc. Japan.} \textbf{8}(1956), 156--166.

\bibitem{Pom64} Ch. Pommerenke,
\textrm{Linear-invariante familien analytischer funktionen. I},
\textit{Math. Ann.} {\bf 155}(1964), 108--154.

\bibitem{Pomm} Ch. Pommerenke,
\textit{Univalent functions},
Vandenhoeck and Ruprecht, G\"ottingen, 1975.

\bibitem{PonRasi2013} S. Ponnusamy and A. Rasila,
Planar harmonic and quasiregular mappings,
Topics in Modern Function Theory (Editors St.Ruschewyeh and S. Ponnusamy):
Chapter in CMFT, \textit{RMS-Lecture Notes Series No.} \textbf{19}, 2013, pp. 267--333.

\bibitem {Schau2000} L. E. Schaubroeck,
\textrm{Subordination of planar harmonic functions},
\textit{Complex Variables} {\bf 41}(2000), 163--178.

\bibitem{Sheil90} T.~Sheil-Small,
\textrm{Constants for planar harmonic mappings},
\textit{J. London Math. Soc.} {\bf 42}(1990), 237--248.

\bibitem{SoStarSzy2011} M. Sobczak-Knec, V. V. Starkov, and J. Szynal,
\textrm{Old and new order of linear invariant family of harmonic mappings and the bound for Jacobian},
\textit{Annales Universitatis Mariae Curie-Sk{\l}odowska. Sectio A. Mathematica,}
\textbf{V. LXV}(2)(2011), 191--202.

\bibitem{Star90}
V.V. Starkov, \textrm{Equivalent definitions of linear invariant families},
\textit{Materialy XI Konf. Szkolenovej z Teorii  Zagadnien Ekstremalnych}, Lodz. 1990, P. 34--38 (in Polish).

\bibitem{Star93} V. V. Starkov,
\textrm{Harmonic locally quasiconformal mappings},
\textit{Trudy of Petrozavodsk State University. Math.} \textbf{1}(1993), 61--69 (in Russian).

\bibitem{Star95} V. V. Starkov,
\textrm{Harmonic locally quasiconformal mappings},
\textit{Annales Universitatis Mariae Curie-Sk{\l}odowska. Sectio A. Mathematica, } \textbf{14}(1995), 183--197.

\bibitem{S2} V. V. Starkov,
Univalence disks of harmonic locally quasiconformal mappings and harmonic Bloch functions,
\textit{Siberian Math. J. } {\bf 38}(1997), 791--800.

\end{thebibliography}
\end{document}